\documentclass[10pt]{article}
\usepackage{amsmath,amssymb}

\newtheorem{theorem}{Theorem}[section]
\newtheorem{lemma}[theorem]{Lemma}
\newtheorem{proposition}[theorem]{Proposition}

\newtheorem{example}[theorem]{Example}
\newtheorem{definition}[theorem]{Definition}
\newtheorem{corollary}[theorem]{Corollary}
\newtheorem{remark}[theorem]{Remark}

\newcommand{\qed}{\enspace\vrule height6pt width4pt depth2pt}
\newenvironment{proof}{\par\noindent{\bf Proof.}}{$\qed$\par\bigskip}

\newcommand{\Z}{{\mathbb Z}}

\newcommand{\R}{{\mathbb R}}

\newcommand{\Fa}{\operatorname{Fa}}

\newcommand{\FaM}{\operatorname{FaM}}
\newcommand{\Sym}{\operatorname{Sym}}
\newcommand{\gr}{\operatorname{gr}}
\newcommand{\Aut}{\operatorname{Aut}}
\newcommand{\Spec}{\operatorname{Spec}}
\newcommand{\hth}{\operatorname{ht}}

\date{}
\begin{document}

\title{Monoids of IG-type and Maximal Orders}
\author{ Isabel
Goffa\footnote{Research funded by a Ph.D grant of
the Institute for the Promotion of Innovation
through Science and Technology in Flanders
(IWT-Vlaanderen).} \and Eric
Jespers\footnote{Research partially supported by the
Onderzoeksraad of Vrije Universiteit Brussel, Fonds
voor Wetenschappelijk Onderzoek (Flanders) and
Flemish-Polish bilateral agreement BIL2005/VUB/06.}
}

\maketitle

\begin{abstract}
Let $G$ be a finite group that acts on an abelian
monoid $A$. If  $\phi :A \rightarrow G$ is a map
so that $\phi (a\phi (a)(b)) =\phi (a)\phi(b)$,
for all $a,b\in A$, then the submonoid $S=\{
(a,\phi (a))\mid a\in A\}$ of the associated
semidirect product $A\rtimes G$ is said to be a
monoid of $IG$-type. If $A$ is a finitely
generated free abelian monoid of rank $n$  and
$G$ is a subgroup of the symmetric group
$\Sym_{n}$ of degree $n$, then these monoids
first appeared in the work of Gateva-Ivanova and
Van den Bergh (they are called monoids of I-type)
and later in the work of Jespers and
Okni\'{n}ski. It turns out that their associated
semigroup algebras share many properties with
polynomial algebras in finitely many commuting
variables.

In this paper we first note that finitely
generated monoids $S$ of IG-type are epimorphic
images of monoids of I-type and their algebras
$K[S]$  are Noetherian and satisfy a polynomial
identity.  In case the group of fractions
$SS^{-1}$ of $S$ also is torsion-free then it is
characterized when $K[S]$ also is  a maximal
order. It turns out that they often are, and
hence these algebras again share arithmetical
properties with natural classes of commutative
algebras. The characterization is in terms of
prime ideals of $S$, in particular $G$-orbits of
minimal prime ideals in $A$ play a crucial role.
Hence, we first describe the prime ideals of $S$.
It also is described when the group $SS^{-1}$ is
torsion-free.
\end{abstract}

\section{Introduction}

In \cite{gateva} Gateva-Ivanova and Van den Bergh
introduced a new class of monoids $T$, called
monoids of I-type, with the aim of constructing
non-commutative algebras that share many
properties with  polynomial algebras in finitely
many commuting variables. In particular, the
semigroup algebras $K[T]$  are Noetherian maximal
orders that satisfy a polynomial identity.
Moreover, these monoids are intimately connected
with set theoretic solutions of the quantum
Yang-Baxter equation and Bieberbach groups. In
this paper, we consider a much wider class of
semigroups $S$ and show that often their algebras
$K[S]$  still are Noetherian maximal orders that
satisfies a polynomial identity.  Earlier recent
results on the construction of such algebras can
be found in \cite{jes,okni,jeswang}, as well as
an extensive literature on the topic.

To put things into context, we first recall the
definition of a monoid of I-type. By $\FaM_{n}$
we denote the free abelian monoid of rank $n$
with basis $\{u_{1},\ldots ,u_{n}\}$. A monoid
$S$, generated by a set $X = \{x_{1},\ldots
,x_{n}\}$, is said to be of left I-type if there
exists a bijection (called a left I-structure)
$v: \FaM_{n} \rightarrow S$ such that $v(1) = 1$
and $\{v(u_{1}a),\ldots ,v(u_{n}a)\} =
\{x_{1}v(a),\ldots ,x_{n}v(a)\}$, for all $a \in
\FaM_{n}$. Similarly one defines monoids of right
I-type. In \cite{gateva} it was shown that a
monoid $S$ of left I-type has a presentation $S =
\langle x_{1},\ldots ,x_{n}\mid R\rangle$, where
$R$ is a set of $\left(\begin{array}{c} n\\ 2\\
\end{array}\right)$ defining relations of the type
$x_{i}x_{j} = x_{k}x_{l}$, so that every word
$x_{i}x_{j}$ with $1 \leq i,j \leq n$ appears at
most once in one of the relations. Hence, one
obtains an associated bijective map $r: X \times
X \rightarrow X \times X$, defined by
$r(x_{i},x_{j}) = (x_{k},x_{l})$ if $x_{i}x_{j} =
x_{k}x_{l}$ is a defining relation for $S$,
otherwise one defines $r(x_{i},x_{j}) =
(x_{i},x_{j})$. For every $x\in X$, denote by
$f_{x}: X \rightarrow X$ and by $g_{x}: X
\rightarrow X$ the mappings defined by
$f_{x}(x_{i}) = p_{1}(r(x,x_{i}))$ and
$g_{x}(x_{i}) = p_{2}(r(x_{i},x))$, where $p_{1}$
and $p_{2}$ denote the projections onto the first
and second component respectively. So,
$r(x_{i},x_{j}) =
(f_{x_{i}}(x_{j}),g_{x_{j}}(x_{i}))$. One says
that $r$ (or simply $S$) is left non-degenerate
if each $g_{x}$ is bijective. In case each
$f_{x}$ is bijective then $r$ (or $S$) is said to
be right non-degenerate. Also, one says that $r$
is a set theoretic solution of the Yang-Baxter
equation if $r_{1}r_{2}r_{1} = r_{2}r_{1}r_{2}$,
where $r_{i}: X^{m}\rightarrow X^{m}$ is defined
as $id_{X^{i-1}} \times r \times id_{X^{m-i-1}}$
and $id_{X^{j}}$ denotes the identity map on the
Cartesian product $X^{j}$.

In \cite{gateva} the equivalence of the first two
statements of the following theorem has been
proven. The equivalence with the third statement
has been proven in \cite{jes}.

\begin{theorem}\label{elem}
The following conditions are equivalent for a
monoid $S$.
\begin{enumerate}
\item $S$ is a monoid of left I-type.
\item
$S$ is finitely generated, say by $x_{1},\ldots
,x_{n}$, and is defined by
$\left(\begin{array}{c} n\\ 2\\
\end{array}\right)$ homogeneous relations of the
form $x_{i}x_{j} = x_{k}x_{l}$ so that every word
$x_{i}x_{j}$ with $1 \leq i,j \leq n$ appears at
most once in one of the relations and the
associated bijective map $r$ is a solution of the
Yang-Baxter equation and is left non-degenerate.
\item
$S$ is a submonoid of a semi-direct product of a
free abelian monoid $\FaM_{n}$ of rank $n$ and a
symmetric group of degree $n$, so that the
projection onto the first component is bijective.
That is, $S = \{(a,\phi(a))\mid a\in \FaM_{n} \}$
where $\phi$ is a mapping from $\FaM_{n}$ to
$\Sym_{n}$ so that
\begin{eqnarray}\label{2}\phi(a)\phi(b) =
\phi(a\phi(a)(b)),\end{eqnarray} or equivalently
\begin{eqnarray}\label{2'} \phi(ac) =
\phi(a) \phi (\phi(a)^{-1}(c)),\end{eqnarray} for
all $a,b,c \in \FaM_{n}$.
\end{enumerate}
\end{theorem}
It follows that a monoid is of left I-type if and
only if it is of right I-type. Such monoids are
simply called monoids of I-type (as in
\cite{jes}).

Note that the above mentioned semi-direct product
$ \FaM_{n}\rtimes \Sym_{n}$ is defined via the
natural action of $\Sym_{n}$ on a chosen basis
$\{u_{1},\ldots ,u_{n}\}$ of the free abelian
monoid $\FaM_{n}$, that is, $\phi(a)(u_{i}) =
u_{\phi(a)(i)}$. Let $\Fa_{n}$ denote the free
abelian group with the same basis. Then, the
monoid $S$ has a group of quotients $SS^{-1}$
contained in $\Fa_{n}\rtimes \Sym_{n}$ and
$SS^{-1} = \{(a,\overline{\phi}(a))\mid a\in
\Fa_{n}\}$, where $\overline{\phi}:
\Fa_{n}\rightarrow \Sym_{n}$ is a mapping that
extends the map $\FaM_{n}\rightarrow \Sym_{n}$
and it also satisfies (\ref{2}). In \cite{jes}
such groups are called groups of I-type. In
\cite{gateva} and \cite{jes} it is shown that
$SS^{-1}$   is a solvable Bieberbach group, that
is, $SS^{-1}$ is a finitely generated solvable
torsion-free group. These groups also have been
investigated by Etingof, Guralnick, Schedler and
Soloviev in \cite{eting-gur,eting}, where they
are called structural groups.

Gateva-Ivanova and Van den Bergh \cite{gateva},
proved that the semigroup algebra of such a
monoid shares a lot of properties with
commutative polynomial algebras in finitely many
variables. In particular, it is a Noetherian
domain that satisfies a polynomial identity and
it is a maximal order.

Jespers and Okni\'{n}ski in \cite{okni}
investigated when an arbitrary semigroup algebra
satisfies these latter properties. The
assumptions on $S$ say that $K[S]$ is a
Noetherian domain that satisfies a polynomial
identity. For details we refer to \cite{okni} and
\cite{boekok} (see also the introduction of
section 3).

\begin{theorem}\label{tool}
Let $K$ be a field and $S$ a submonoid of a
torsion-free finitely generated abelian-by-finite
group. The monoid algebra $K[S]$ is a Noetherian
maximal order if and only if the following
conditions are satisfied:
\begin{enumerate}
\item $S$ satisfies the ascending chain
condition on one sided ideals,
\item
$S$ is a maximal order in its group of quotients
$H=SS^{-1}$,
\item
for every minimal prime $P$ in $S$, $$S_{P} =
\{g\in H\mid Cg \subseteq S \quad \text{for some
$H$-conjugacy class C of $H$}$$ $$\text{contained
in S and with $C \nsubseteq P$}\}$$ has only one
minimal prime ideal.
\end{enumerate}
\end{theorem}

It is worth mentioning that Brown in \cite{brown}
proved that, for a field $K$, a group algebra
$K[G]$ of a torsion-free polycyclic-by-finite
group $G$ is a maximal order in its classical
ring of quotients (which is a domain). A
characterization of commutative semigroup
algebras $K[A]$ that are Noetherian domains and
maximal orders can be found in \cite{gilmer}. It
turns out that $K[A]$ is such an algebra if and
only if $A$ is finitely generated and a maximal
order in its torsion-free group of quotients
$AA^{-1}$ (see also the comments given in Section
3). Extensions of this result to Krull orders
have been proved by Chouinard \cite{choui}.
\section{Monoids of IG-type}
We begin with introducing the  larger class of
monoids of interest. Let $G$ be a group and $A$ a
monoid. Recall that $G$ is said to act on $A$ if
there exists a monoid morphism $\varphi :
G\rightarrow \Aut (A)$. The associated
semi-direct product $A\rtimes_{\varphi} G$ we
often simply denote by $A\rtimes G$.

\begin{definition}
Suppose $G$ is a finite group acting on a
cancellative abelian monoid $A$. A submonoid $S$
of $A\rtimes G$ so that the natural projection on
the first component is bijective is said to be a
monoid of IG-type. Thus,
 $$S = \{(a,\phi(a))\mid a\in A \},$$ with
$\phi : A\rightarrow G$  a mapping satisfying
(\ref{2}). (We denote the action of $g\in G$ on
$a\in A$ as $g(a)$.)
\end{definition}

 Note
that for every $a\in A$,
$$(a,\phi(a))^{\left|G\right|} = (a\phi (a)\cdots
\phi (a)^{\left|G\right| - 1} (a), \phi
(a)^{\left|G\right|}) = (a\phi (a)\cdots \phi
(a)^{\left|G\right| - 1} (a),1).$$ It follows
that $b=a\phi(a)\cdots \phi(a)^{\left|G\right| -
1}(a)\in A$ is such that $\phi(b) = 1$ and, for
every $a^{-1}b_{1}\in AA^{-1}$ (the group of
quotients of $A$), we have $a^{-1}b_{1} =
b^{-1}(\phi(a)\cdots \phi(a)^{\left|G\right|
-1}(a)b_{1})$. So, any element of $AA^{-1}$ can
be written as $a^{-1}b$ with $a,b\in A$ and
$\phi(a) = 1$. Furthermore, if $a_{1}^{-1}b_{1} =
a_{2}^{-1}b_{2}$, with $a_{i}, b_{i}\in A$ and
$\phi(a_{i}) = 1$, then by equation (\ref{2}), it
is easily verified that $\phi(b_{1}) =
\phi(b_{2})$. We hence can extend the action of
$G$ onto $A$ to an action of $G$ onto $AA^{-1}$
and thus obtain a mapping $\overline{\phi} :
AA^{-1}\rightarrow G$ so that
$$\overline{\phi}(a^{-1}b) = \phi(b),$$ for every
$a,b\in A$ with $\phi(a) = 1$. This mapping again
satisfies (\ref{2}).

Hence $S$ is a submonoid of the group
$AA^{-1}\rtimes G.$ Since $AA^{-1}\rtimes G$ is
abelian-by-finite and $S$ is cancellative, Lemma
7.1 in \cite{boekok} yields that $S$ has a group
of fractions $SS^{-1}\subseteq AA^{-1}\rtimes G$.
Furthermore, because of Theorem 15 in
\cite{boekok}, the algebra  $K[S]$ satisfies a
polynomial identity, and if $K[S]$ is prime then
$SS^{-1} = SZ(S)^{-1}$, with $Z(S)$ the center of
$S$ \cite{jeswang}.

We claim that the natural projection of
$SS^{-1}\rightarrow AA^{-1}$ is a one-to-one
mapping. Indeed, if $a_{1}^{-1}b_{1} =
a_{2}^{-1}b_{2}$ (with $a_{i},b_{i}\in A,
\phi(a_{i}) = 1$) then, by (\ref{2}),
$\phi(b_{1}) = \phi(b_{2})$. So,
$$(a_{1}^{-1}b_{1},\overline{\phi}(a_{1}^{-1}b_{1}))
= (a_{1}^{-1}b_{1},\phi(b_{1})) =
(a_{2}^{-1}b_{2}, \phi(b_{2})) =
(a_{2}^{-1}b_{2},
\overline{\phi}(a_{2}^{-1}b_{2})).$$ This proves
the injectiveness. The surjectiveness follows
from the fact that
\begin{eqnarray*}
(a_{1}^{-1}b_{1},
\overline{\phi}(a_{1}^{-1}b_{1})) &=&
(a_{1}^{-1}b_{1},\phi(b_{1}))\\
 &=&
(a_{1},1)^{-1}(b_{1},\phi(b_{1}))\\
 & =&
(a_{1},\phi(a_{1}))^{-1}(b_{1},\phi(b_{1}))\in
SS^{-1}, \end{eqnarray*} for $a_{i},b_{i}\in A$
with $\phi(a_{i}) = 1$. Groups of the type
$SS^{-1}$ we call groups of IG-type. So we have
shown the following.

\begin{corollary}
A group $H$ is of IG-type if and only if $H$ is a
subgroup of a semi-direct product $A\rtimes G$ of
a finite group $G$ with an abelian group $A$ so
that
 $$H = \{(a,\overline{\phi}(a))\mid a\in A\}$$
and
 $$\overline{\phi}(a\overline{\phi}(a)(b)) =
\overline{\phi}(a)\overline{\phi}(b),$$ for all
$a,b\in A$. Of course, such a group is
abelian-by-finite.
\end{corollary}

A subset $B$ of $A$ is said to be
$\phi$-invariant if $\phi (a)(B)=B$ for all $a\in
A$. In case $B$ is a subgroup of $A$ then this
condition is equivalent with $B$ being a normal
subgroup of $SS^{-1}$.

Note also that if $S=\{ (a,\phi (a) \mid a\in A\}
\subseteq A\rtimes G$ is a monoid of IG-type,
with $A$ an abelian monoid and $G=\{ \phi (a)
\mid a\in A\}$ a finite group then $\prod_{\phi
(a)\in G} \phi (a) (b)$ is an invariant element
of $A$, for every $b\in A$. (In this case we will
also use $G$-invariant as $\phi$-invariant) It
follows that every element of $SS^{-1}$ can be
written as $(z,1)^{-1}(a,\phi (a))$ with $z,a\in
A$ and $z$ and invariant element in $A$. So
$(z,1)$ is a central element of $S$.

We now describe when the semigroup algebra $K[S]$
of a monoid of IG-type is Noetherian. This easily
can be deduced from the following Lemma and the
recent result of Jespers and Okni\'{n}ski proved
in \cite{noet} which says that, for a submonoid
$T$ of a polycyclic-by-finite group, the
semigroup algebra $K[T]$ is left Noetherian if
and only if $K[T]$ is right Noetherian, or
equivalently, $T$ satisfies the ascending chain
condition on left ( or right) ideals. An algebra
which is left and right Noetherian we simply call
Noetherian. However, for completeness' sake we
include a simple proof for the monoids under
consideration. The subgroup generated by a set
$X$ of elements in a group  $G$ is denoted $\gr
(X)$. By  $\langle X \rangle$ we denote  the
monoid generated by $X$.

\begin{lemma}\label{laatste}
Let $A = \langle u_{1},\ldots ,u_{n}\rangle$ be a
finitely generated abelian monoid, $G$ a finite
group acting on $A$. Let $S = \{(a,\phi(a))\mid
a\in A\} \subseteq A\rtimes G$ be a monoid of
IG-type. Put   $B = \{ \phi (a)(u_{i}) \mid a\in
A,\; 1\leq i \leq n\}$. Then the following
conditions hold:
\begin{enumerate}
\item $G$ acts on the set $B$, that is, $\phi (a) (B)=B$,
for all $a\in A$.
\item $S = \langle (b,\phi (b)) \mid b\in B \rangle$.
\item for some divisor $k$ of $|G|$, the subgroup
$\gr \{ (b^{k},1) \mid b\in B\} $ is normal and
of finite index in $SS^{-1}$.
\item
$S=\bigcup_{f\in F} \langle (b^{k},1)\mid b\in B
\rangle (f,\phi (f))$ and $(f,\phi (f) \langle
(b^{k},1) \mid b\in B \rangle = \langle
(b^{k},1)\mid b\in B \rangle (f,\phi(f))$, for
some finite subset $F$ of $A$.
\end{enumerate}
\end{lemma}
\begin{proof}
The first and second part  follow at once from
the equalities (\ref{2}) and (\ref{2'}).   Put $N
= \{a\in AA^{-1}\mid \overline{\phi}(a) = 1\}$,
the kernel of the natural homomorphism
$SS^{-1}\rightarrow G$. So, $N$ is an abelian
subgroup of finite index $k$ in $AA^{-1}$, with
$k$ a divisor of $|G|$. It follows that
$\phi(a^{k}) = 1$, for any $a\in A$ and that  the
abelian monoid $C= \langle (b^{k},1) \mid b\in
B\rangle$ is contained in $S$ and its group of
quotients $CC^{-1}= \gr \{ (b^{k},1) \mid b\in B
\}$ is normal and of finite index in $SS^{-1}$.
This proves the third part. Part four is now also
clear.
\end{proof}

If, in the previous lemma,  $U(A) = \{1\}$ then
one can take $\{ u_{1},\ldots , u_{n}\}$ to be
the set of indecomposable elements, that is, the
set consisting of those elements $f\in A$ so that
$Af$ is a maximal principal ideal. Indeed, since
$A$ is finitely generated, we know that  $A$
satisfies the ascending chain condition on
ideals. Hence, $A$ has finitely many
indecomposable elements, say $u_{1},\ldots ,
u_{n}$,  and $A = \langle u_{1},\ldots
,u_{n}\rangle$ (\cite{pbf}). Clearly any
automorphism of $A$ permutes the indecomposable
elements. It follows that $S=\langle (u_{1},\phi
(u_{1})),\ldots , (u_{n},\phi (u_{n}))\rangle$.

Let $S$ be a submonoid of a group $G$. Assume $N$
is a normal subgroup of $G$. If $N\subseteq S$
then, as in group theory we denote by $S/N$ the
monoid consisting of the cosets $sN=Ns$, with
$s\in S$.

\begin{proposition} \label{noeth}
Let $A$ be an abelian monoid and $G$ a finite
group acting on $A$. Let $S = \{(a,\phi(a))\mid
a\in A\} \subseteq A\rtimes G$ be a monoid of
IG-type. Then, the semigroup algebra $K[S]$ is
Noetherian if and only if the abelian monoid $A$
is finitely generated, or, equivalently, $S$ is
finitely generated.
\end{proposition}

\begin{proof}
Suppose that $K[S]$ is right Noetherian. Because
$K[L]$ is a right ideal of $K[S]$ for every right
ideal $L$ of $S$, it follows easily that $S$
satisfies the ascending chain condition on right
ideals. Consequently, also $A$ satisfies the
ascending chain on right ideals. Indeed, if $L$
is a right ideal of $A$, then $\lambda_{L} =
\{(a, \phi(a))\in S \mid a \in L\}$ is a right
ideal of $S$ and $\lambda_{L}\subset
\lambda_{L'}$ if and only if $L \subset L'$. So
$A$ is an abelian and cancellative monoid that
satisfies the ascending chain condition on
ideals. Hence, so is the monoid $A/U(A)$. Because
$U(A/U(A))=\{ 1 \}$, it follows (see the remark
above)  that $A/U(A)$ is finitely generated by
its indecomposable elements. Clearly $U(S)=\{
(a,\phi (a))\mid a\in U(A)\}$ and because
$IK[S]\cap K[U(S)]=I$ for any right ideal $I$ of
$K[U(S)]$, it follows that the group algebra
$K[U(S)]$ is Noetherian. Hence it is well known
that $U(S)$ is a finitely generated monoid.
Consequently, $U(A)$ and thus also $A$ is
finitely generated.

For the converse, suppose that $A = \langle
u_{1},\ldots   ,u_{n}\rangle$ is finitely
generated.
>From Lemma~\ref{laatste} it follows that the algebra
$K[S]$ is a finite module over the commutative
Noetherian algebra $K[\langle (b^{k},1) \mid b\in
B\rangle]$. Hence $K[S]$ is Noetherian.
\end{proof}

We now give a link with monoids of I-type by
proving another characterization of finitely
generated monoids $S$ of IG-type.

\begin{theorem}\label{epi}
A finitely generated monoid $S$ is of IG-type if
and only if there exists a monoid of $I$-type
$T=\{ (x,\psi (x)) \mid x\in \FaM_{m}\}\subseteq
\FaM_{m}\rtimes \Sym_{m}$ and a subgroup $B$ of
$\Fa_{m}$ that is $\psi$-invariant so that
$S\cong TB/B$.
\end{theorem}

\begin{proof}
Assume $S=\{ (a,\phi (a) ) \mid a\in A\}\subseteq
A\rtimes G$ is a finitely generated monoid of
IG-type, where $G$ is a finite group acting on
the finitely generated abelian monoid $A=\langle
u_{1}, \ldots , u_{n}\rangle $. Of course we may
assume that $G=\{ \phi (a) \mid a\in A\}$.

Let $m=n|G|$ and let $\FaM_{m}$ be the free
abelian monoid of rank $m$ with basis the set
$M=\{ v_{g,i} \mid g\in G, \; 1\leq i \leq n\}$.
Clearly the mapping $f: \FaM_{m}\rightarrow A$
defined by $f(v_{g,i})=gu_{i}$ is a monoid
epimorphism. For $x\in \FaM_{m}$ define a mapping
$\psi (x):M \rightarrow M$ by $\psi (x)
(v_{g,i})= v_{\phi (f(x)) g,i}$. Then $\psi (x)
\in \Sym_{m}$, $f(\psi (x) (y)) = \phi (f(x))
(f(y))$ and  $\psi (x\psi (x)(y)) =\psi (x) \psi
(y)$, for any $x,y\in \FaM_{m}$. So $T=\{ (x,\psi
(x)) \mid x\in \FaM_{m}\}$ is a monoid of
$I$-type contained in $\FaM_{m}\rtimes \Sym_{m}$.
Furthermore, $f^{e}:T\rightarrow S$ defined by
$f^{e}((x,\psi (x))) =(f(x),\phi (f(x)))$ is a
monoid epimorphism. Its extension to an
epimorphism $TT^{-1}\rightarrow SS^{-1}$ we also
denote by $f^{e}$. Let $B=\ker (f^{e})$. Clearly,
if  $(x,\psi (x))\in B$ then  $\phi (f(x))=1$ and
thus $\psi (x)=1$. Thus, $B\subseteq \Fa_{m}$ and
$B$ is $\psi$-invariant. So $TB=BT$ is a
submonoid of $TT^{-1}$ and $TB/B \cong S$. This
proves the necessity of the conditions.

Conversely, assume  $T=\{ (x,\psi (x)) \mid x\in
\FaM_{m}\}\subseteq \FaM_{m}\rtimes \Sym_{m}$ is
a monoid of $I$ type and  $B$ is a subgroup  of
$\Fa_{m}$ that is $\psi$-invariant. Let
$A=\FaM_{m}B/B$ and let $f:\FaM_{m}\rightarrow A$
be the natural monoid epimorphism. Because of
(\ref{2'}) we get that $\psi (x)=\psi (y)$ if
$f(x)=f(y)$. Hence, for each $a=f(x)$ the mapping
$\phi (a):A\rightarrow A$ given by $\phi
(a)(f(y)) =f(\psi (x)(y))$ is a well defined
bijection of finite order. Furthermore, $\phi
(a\phi (a)(b))=\phi (a)\phi (b)$ for all $a,b\in
A$. Hence $S=\{ (a,\phi (a)) \mid a\in A\}$ is a
monoid of IG-type contained in $A\rtimes G$,
where $G=\{ \phi (a)\mid a\in A\}$. The mapping
$TB\rightarrow S$ defined by mapping $(x,\psi
(x))$ onto $(a,\phi (a))$ is a monoid epimorphism
and it easily follows that this map induces an
isomorphism between $TB/B$ and $S$.
\end{proof}

We note that the proposition can be formulated
using congruence relations as follows. A finitely
generated monoid $S$ is of IG-type if and only if
there exists a monoid of $I$-type  $T =
\{(a,\psi(a))\mid a \in \FaM_{m}\} \subseteq
\FaM_{m}\rtimes \Sym_{m}$  and there exists a
congruence relation $\rho$ on $\FaM_{m}$ with
\begin{equation} \label{geg}
a\;  \rho \; b  \quad \mbox{ implies } \quad
\psi(a) = \psi(b)\mbox{ and } \psi(x)(a) \; \rho
\; \psi(x)(b),
\end{equation}
for every $a,b,x \in \FaM_{m}$, and so that
$S\cong T/ \overline{\rho}$  where
$\overline{\rho}$ is the congruence relation on
$T$ defined by $(a,\psi(a))\; \overline{\rho}\;
(b,\psi(b))$ if and only if $a\; \rho \; b$, for
$a,b \in \FaM_{m}$.

We remark that many monoids of IG-type are not of
$I$-type. Indeed, suppose $S=\{ (a,\phi (a))\mid
a\in A\}$ is a monoid of IG-type with $A$ a
finitely generated monoid so that $U(A)=\{ 1\}$.
Let $\{ u_{1}, \ldots , u_{n}\}$ be the set of
indecomposable elements of $A$. So $A=\langle
u_{1}, \ldots , u_{n}\rangle$. It follows that
the elements $(u_{i},\phi(u_{i}))$ are the unique
indecomposable elements of $S$, that is, they can
not be decomposed as a product of two
non-invertible elements. So $S$ also has $n$
indecomposables. In particular, the number of
indecomposables in a monoid $T$ of $I$-type
equals the torsion-free rank of any abelian
subgroup of finite index in $TT^{-1}$. So, if the
torsion-free rank of $AA^{-1}$ is strictly
smaller than $n$ then $S$ is not of $I$-type.

\section{Torsion-freeness of Groups of IG-type}

We recall some notation and terminology on
maximal orders (see for example \cite{okni}). A
cancellative monoid $S$ which has a left and
right group of quotients $G$ is called an order.
Such a monoid $S$ is called a maximal order if
there does not exist a submonoid $S'$ of $G$
properly containing $S$ and such that
$aS'b\subseteq S$ for some $a,b \in S$. For
subsets $A,B$ of $G$ put $(A:_{l} B) = \{g\in
G\mid gB\subseteq A\}$ and $(A:_{r} B) = \{g\in
G\mid Bg\subseteq A\}$. It turns out that $S$ is
a maximal order if and only if $(I :_{l} I) = (I
:_{r} I) = S$ for every fractional ideal $I$ of
$S$. The latter means that $SIS \subseteq I$ and
$cI, Id \subseteq S$ for some $c,d \in S$. If $S$
is a maximal order, then $(S:_{l}I) = (S:_{r}I)$
for any fractional ideal $I$. One simply denotes
this fractional ideal by $(S : I)$ or by
$I^{-1}$. Recall that $I$ is said to be
divisorial if $I = I^{*}$, where $I^{*} = (S : (S
: I))$. The divisorial product $I * J$ of two
divisorial ideals $I$ and $J$ is defined as
$(IJ)^{*}$.

Also recall that a cancellative monoid $S$ is
said to be a Krull order if and only if $S$ is a
maximal order satisfying the ascending chain
condition on integral divisorial ideals, that is,
fractional ideals contained in $S$. In this case
the set $D(S)$ of divisorial ideals is a free
abelian group for the
* operation. If $G$ is abelian-by-finite, then every
ideal of $S$ contains a central element. In this
case, it follows that the minimal primes of $S$
form a free basis for $D(S)$. The positive cone
of this group (with respect to this basis) is
denoted by $D(S)^{+}$.

In this section, we investigate periodic elements
of a monoid $S = \{(a,\phi(a))\mid a\in A\}$ of
IG-type, and we will restrict our attention to
the case that $A$ is a finitely generated maximal
order (and hence a Krull order) with trivial unit
group, $AA^{-1}$ is torsion-free and the action
of $G=\{ \phi (a)\mid a\in A\}$ on $A$ is
faithful. Because of the latter condition we may
consider $G$ as a subgroup of the automorphism
group of $A$. Since $A$ is finitely generated, we
know that $A$ has only finitely many minimal
prime ideals and every prime ideal is a union of
minimal prime ideals. Recall from Theorem~37.5 in
\cite{pass} that, as $SS^{-1}$ is
polycyclic-by-finite, $K[S]$ (or equivalently
$K[SS^{-1}]$) is a domain if and only if
$SS^{-1}$ is torsion-free.

\begin{proposition}\label{ok}
Let $A$ be an abelian cancellative monoid. Assume
that $A$ is a finitely generated maximal order
and $U(A) = \{1\}$. If $S = \{(a,\phi(a))\mid a
\in A\} \subseteq A \rtimes G$ is a monoid of
IG-type and the action of $G=\{ \phi (a) \mid
a\in A\}$ is faithful, then $$S \cong
(D(A)^{+}\rtimes G)\cap SS^{-1},$$ where $D(A)$
is the divisor class group of $A$ and $G$ is a
subgroup of the permutation group of the minimal
primes of $A$.
\end{proposition}
\begin{proof}
As $A$ is a Krull order, we know that $D(A)$ is a
free abelian group with the set $\Spec^{0}(A)$
consisting of the  minimal primes of $A$ as a
free generating set. Of course, for each $a\in
A$, $\phi(a)$ induces an automorphism on
$\Spec^{0}(A)$, and thus also on $D(A)$. We
denote this again by $\phi(a)$. It follows that,
if $I$ is an ideal of $S$, then $\phi(a)(I^{*}) =
(\phi(a)(I))^{*}$. We thus obtain a morphism $G
\rightarrow \Sym (\Spec^{0}(A))$. This mapping is
injective. Indeed, suppose  $\phi(a)$ is the
identity map on $\Spec^{0}(A)$, with $a\in A$.
For $c\in A$ the ideal $Ac$ is divisorial. hence
$A\phi(a)(c)= \phi(a) (Ac)=Ac$. Since, by
assumption $U(A)=\{1 \}$ it follows that $\phi
(a)(c)=c$. Hence it follows that $\phi (a)=1$.
Because, also by assumption, the action of $G$ on
$A$ is faithful, it follows that $\phi (a)$ is
the identity map on $G$, as desired.

Again, because $U(A)=\{ 1\}$, we get a monoid
morphism $$S \rightarrow D(A)^{+}\rtimes G:
(a,\phi(a))\mapsto (aA,\phi(a)).$$ So,
identifying $S$ with its image  in
$D(A)^{+}\rtimes G$ (and also $SS^{-1}$ with its
image in
$\{(a^{-1}bA,\overline{\phi}(a^{-1}b))\mid a,b\in
A\} \subseteq D(A)\rtimes G$), we get that $$S
\subseteq (D(A)^{+}\rtimes G)\cap SS^{-1}.$$
Conversely, if $b\in AA^{-1}$ and $(bA,
\overline{\phi}(b))\in D(A)^{+}\rtimes G$, then
$bA \subseteq A$ and thus $b\in A$. Hence $(bA,
\overline{\phi}(b))\in S$.
\end{proof}

Note that this characterization is a
non-commutative version of result of Chouinard
that describes commutative cancellative
semigroups that are Noetherian maximal orders or
more generally Krull orders \cite{choui}.

To investigate the torsion-freeness, we need the
following Theorem. The authors would like to
thank Karel Dekimpe for the proof of this result
\cite{karel}.
\begin{theorem}\label{karel}
Let $H$ be a group of affine transformations such
that $H\cap \R^{n}$ (the subgroup of pure
translations) is of finite index in $H$. Then the
following properties are equivalent.
\begin{enumerate}
\item $H$ is torsion-free.
\item
The action of $H$ on $\R^{n}$ is fixed-point
free, that is, if $g \cdot a = a$ for some $a\in
\R^{n}$ and $g\in H$, then $g = 1$.
\end{enumerate}
\end{theorem}
\begin{proof}
Suppose that the action of $H$ on $\R^{n}$ has a
fixed-point. Let therefore $h\neq 1$ and $x\in
R^{n}$ be such that $h \cdot x = x$. Then, also
$h^{k}\cdot x = x$ for every $k \in \Z$. Because
$H\cap \R^{n}$ is of finite index in $H$ there
exists a $k \rangle 0$ such that $h^{k}$ is a
pure translation. But as this translation has a
fixed-point it follows that $h^{k}$ should be
trivial. Therefore, $H$ has torsion.

Conversely, suppose that $H$ is a finite subgroup
of the affine transformations in dimension $n$.
So every element $h$ of $H$ is of the form
$(t_{h}, M_{h})$, with $t_{h}$ the translation
part and $M_{h}$ the linear part. Take $h_{1},
h_{2}\in H$ then: $$(t_{h_{1}h_{2}},
M_{h_{1}h_{2}}) = (t_{h_{1}} +
M_{h_{1}}t_{h_{2}}, M_{h_{1}}M_{h_{2}}).$$
Therefore the map $$\phi: H\rightarrow GL(n,\R):
h \mapsto M_{h},$$ is a group morphism. Hence
$\R^{n}$ is an $H$-module. With this module
structure the map $t:H \rightarrow \R^{n}$
becomes a $1$-cocycle. Because $H$ is finite, we
have that $H^{1}(H, \R^{n}) = 0$ and therefore
the map $t$ is a $1$-coboundary \cite{bab}.
Consequently there exists a $x\in \R^{n}$ with
$t_{h} = x - M_{h}x$ for every $h \in H$ and
therefore we have a fixed-point.
\end{proof}

Note that, if $SS^{-1}$ is torsion-free, then so
is necessarily $AA^{-1}$. Indeed, If $a^{m} = 1$
in $AA^{-1}$, then, by Lemma~\ref{laatste},
$(a\phi(a)(a)\phi(a)^{2}(a)...\phi(a)^{k-1}(a),1)^{m}
= 1 \in SS^{-1}$, for some divisor $k$ of $|G|$.

Assume now that $S = \{(a,\phi(a))\mid a \in A\}$
is a monoid of IG-type with faithful action of
$G=\{ \phi (a) \mid a\in A\}$ on $A$. Suppose
that $AA^{-1}$ is a torsion-free finitely
generated abelian group. So, $SS^{-1}\subseteq
\Z^{k}\rtimes G$ and $G\subseteq Aut(\Z^{k})
\cong GL_{k}(\Z)$. Hence, every element of $G$
can be seen as a $k \times k$-matrix with values
in $\Z$ and the action of $SS^{-1}$ on $\Z^{k}$
can be extended to $\R^{k}$ and can be written
as: $$(a,\phi(a))\cdot b = \phi(a)b + a,$$ where
$a\in \Z^{k}, \phi(a)\in GL_{k}(\Z)$, $b\in
\Z^{k}$ (or $\R^{k}$) and $\phi(a)b$ is given by
the classical matrix multiplication. For
convenience sake we  use the additive notation on
$\Z^{k}$ and $\R^{k}$ (instead of the
multiplicative on $AA^{-1}$). So,
$SS^{-1}\subseteq \R^{k}\rtimes GL_{k}(\R)$.
Thus, $SS^{-1}$ is a group of affine
transformations and every element of $SS^{-1}$ is
of the form $(a,\phi(a))$, where $a$ is the
translation part and $\phi(a)$ the linear part.
Clearly, $(a,A)(b,B) = (a+Ab,AB)$.  As $G$ is a
finite group, we also have that the subgroup of
pure translations is of finite index in
$SS^{-1}$.

As an immediate consequence of Theorem
\ref{karel}, we get that the quotient group
$SS^{-1}\subseteq AA^{-1}\rtimes G$ is
torsion-free if and only if the action of
$SS^{-1}$ on $\R^{n}$ is fixed-point free.

If also $U(A)=\{ 1 \}$ and $A$  is a maximal
order in its free abelian group of quotients
then, by Theorem \ref{ok}, we can extend the
action of $SS^{-1}$ to an action of the
semi-direct product $D(A)\rtimes G$ and $G$ acts
as the symmetric group on the set $\Spec^{0}(A)=
\{P_{1},...,P_{l}\}$. Hence the action can also
naturally be extended to an action of the
semi-direct product $\R^{l}\rtimes G$.

\begin{theorem}\label{torsion-free}
Let $A$ be an abelian cancellative monoid. Assume
$A$ is a finitely generated maximal order with
$U(A) = \{1\}$ and $AA^{-1}$ is torsion-free. If
$S = \{(a,\phi(a))\mid a \in A\} \subseteq A
\rtimes G$ is a monoid of IG-type with faithful
action of $G=\{ \phi (a) \mid a\in A\}$ on $A$,
then an element $(a,\overline{\phi}(a))$ of
$SS^{-1}$ is periodic if and only if there exists
a divisorial ideal $I$ of $A$ such that
$a\overline{\phi}(a)(I) = I$.
\end{theorem}

\begin{proof}
Suppose $(a,\overline{\phi}(a))$ is a periodic
element. So, because of Proposition \ref{ok} and
the proof of Theorem \ref{karel},
$(aA,\overline{\phi}(a))\in D(A)\rtimes \Sym_{l}$
has a fixed point $b$ in $\R^{l}$. Write $aA =
P_{1}^{\alpha_{1}}\ast ...\ast
P_{l}^{\alpha_{l}}$, with $\alpha_{i} \in \Z^{l}$
and $\Spec^{0}(A)= \{P_{1},...,P_{l}\}$. So
$\overline{\phi}(a) \cdot b + \alpha = b$, where
$\alpha = (\alpha_{1},...,\alpha_{l})$. Since
$\overline{\phi}(a)$ acts as a permutation on the
components of $b\in \R^{l}$, it is not so
difficult to see that $\overline{\phi}(a) \cdot
\left\lfloor b \right\rfloor + \alpha =
\left\lfloor b \right\rfloor$, where
$\left\lfloor b \right\rfloor$ is the integral
part of $b$. Hence we have a fixed point in
$\Z^{l}$. This means that
$a\overline{\phi}(a)(I)= I$, where $I =
P_{1}^{\beta_{1}}\ast ...\ast P_{l}^{\beta_{l}}$,
with $(\beta_{1},...,\beta_{l}) = \left\lfloor b
\right\rfloor$.

Conversely, suppose there exists a divisorial
ideal $I$ of $A$ such that
$a\overline{\phi}(a)(I) = I$. It follows that
$$a\; \overline{\phi}(a)(a)\;
\overline{\phi}(a)^{2}(a)\, \cdots \;
\overline{\phi}(a)^{n}(a)\;
\overline{\phi}(a)^{n+1}(I) = I.$$ So if
$\overline{\phi}(a)^{n+1} = 1$, then we obtain
that $a\; \overline{\phi}(a)(a)\; \cdots \;
\overline{\phi}(a)^{n}(a) A = A$. Again because
$U(A)=1$, it follows that  $a\;
\overline{\phi}(a)(a)\; \cdots \;
\overline{\phi}(a)^{n}(a) = 1$. Consequently,
$(a,\overline{\phi}(a))^{n+1} = 1$, as desired.
\end{proof}

We now give concrete examples of monoids of
IG-type that are not of I-type. The first one is
based on an example of an abelian finitely
generated monoid that is considered by Anderson
in \cite{and}.

\begin{example}\label{and}
Let $A = \langle u_{1},u_{2},u_{3},u_{4}\mid
u_{i}u_{j}=u_{j}u_{i}, \; u_{1}u_{2} =
u_{3}u_{4}\rangle $ and let $\left| \right|: A
\rightarrow \Z$ denote the degree function on $A$
defined by $\left|u_{i}\right| = 1$. Put $$\sigma
= \left(\begin{array}{ccc} 0 & 1 & 1 \\ 1 & 0 & 1
\\ 0 & 0 & -1
\end{array}\right) \in Gl_{3}(\Z).$$
The natural action of $\Z_{2} = \langle \sigma
\rangle $ on $\Z^{3} = AA^{-1}=\gr
(u_{1},u_{2},u_{3})$ defines a semi-direct
product $A \rtimes \Z_{2}$. Then
 $$S = \{(a,\phi(a))\mid a\in
A, \;  \phi(a)=1 \text{ if }\; \left|a \right|\in
2\Z, \; \phi(a)=\sigma \text{ if } \left|a
\right|\in 2\Z + 1\}\subseteq A\rtimes \Z_{2}$$
is a monoid of IG-type (which is not of I-type)
and its group of quotients $$ SS^{-1}
=\{(a,\phi(a))\mid a\in \Z^{3},\; \phi(a)=1
\text{ if } \left|a \right|\in 2\Z, \;
\phi(a)=\sigma \text{ if }= \left|a \right|\in
2\Z + 1\}$$ is torsion-free.
\end{example}

\begin{proof}
Because $A = AA^{-1}\cap F^{+}$, the intersection
of the group of quotients  $AA^{-1}$ and the
positive cone of a free abelian group, we know
that $A$ is a maximal order (see
\cite{and,choui}). Clearly $U(A) = \{1\}$ and $A$
has four minimal primes: $Q_{1} = (u_{1},
u_{3})$, $Q_{2} = (u_{1}, u_{4})$, $Q_{3} =
(u_{2}, u_{3})$ and $Q_{4} = (u_{2}, u_{4})$. So,
these minimal primes generate the free abelian
group $D(A)\cong \Z^{4}$. Because of
Theorem~\ref{torsion-free}, to prove that
$SS^{-1}$ is torsion-free we need to show that if
$(a,\phi(a))\in SS^{-1}$ such that $a\phi(a)(I) =
I$ for some divisorial ideal $I$ of $A$ then
$a=1$. Clearly, if $ \left|a \right|\in 2\Z$ then
$\phi(a) = \{1\}$, and thus $aI=I$ implies $a=1$.
So, suppose $a$ has odd degree. Then $a =
u_{1}^{a_{1}}u_{2}^{a_{2}}u_{3}^{a_{3}}$ or $a =
u_{1}^{a_{1}}u_{2}^{a_{2}}u_{4}^{a_{4}}$. We deal
with the former case (the other case is dealt
with similarly). It is readily verified that
$Au_{1} = Q_{1} \ast Q_{2}$, $Au_{2} = Q_{3} \ast
Q_{4}$ and $Au_{3} = Q_{1}\ast Q_{3}$. Write  $I
= Q_{1}^{\gamma_{1}}\ast Q_{2}^{\gamma_{2}}\ast
Q_{3}^{\gamma_{3}}\ast Q_{4}^{\gamma_{4}}$, with
each $\gamma_{i}\in \Z$. Because $\sigma$
interchanges $Q_{1}$ with $Q_{4}$ and $Q_{2}$
with $Q_{3}$, the equality $Aa \ast \phi (a) (I)
= (Aa\phi(a)(I))^{*}=I$ becomes
$$Q_{1}^{a_{1}+a_{3}}\ast Q_{2}^{a_{1}} \ast
Q_{3}^{a_{2}+a_{3}}\ast Q_{4}^{a_{2}}\ast
Q_{4}^{\gamma_{1}}\ast Q_{3}^{\gamma_{2}}\ast
Q_{2}^{\gamma_{3}}\ast Q_{1}^{\gamma_{4}} =
Q_{1}^{\gamma_{1}}\ast Q_{2}^{\gamma_{2}}\ast
Q_{3}^{\gamma_{3}}\ast Q_{4}^{\gamma_{4}}.$$ It
follows that $a_{1}+a_{2}+a_{3}=0$, in
contradiction with the fact that $a$ is of odd
degree.

Note that, as $AA^{-1}$ has torsion-free rank
$3$, while $A$, and therefore also $S$, has $4$
indecomposable elements, it follows from the
remark at the end of Section 2 that $S$ is not of
I-type.
\end{proof}

A second type of examples of monoids of IG-type
that are not of I-type can be constructed as a
natural class of submonoids of a monoid of
I-type. In \cite[Section 4]{jes} an example in
this class is given to show that there exists
monoids $T$ of I-type with a group of fractions
$TT^{-1}$ that is not poly-infinite cyclic.

\begin{example}\label{belvb}
Let $T = \{(a,\phi(a))\mid a\in
\FaM_{n}\}\subseteq \FaM_{n} \rtimes G $ be a
monoid of I-type with $H=\{ \phi (a)\mid a\in
A\}$. Suppose $B$ is a $G$-invariant submonoid of
$\FaM_{n}$. Then, $S = \{(b,\phi(b))\mid b\in
B\}$ is a monoid of IG-type. Note again that if
$B = BB^{-1} \cap \FaM_{n}$ then we know from
\cite{choui} that $B$ is a maximal order.
Clearly, $U(B) = \{1\}$.
\end{example}

We give a concrete example. Let $T = \langle
x_{1},x_{2},x_{3},x_{4}\rangle $ be the monoid
defined by the relations $x_{1}x_{2} =
x_{3}x_{3}$, $x_{2}x_{1} = x_{4}x_{4}$,
$x_{1}x_{3} = x_{2}x_{4}$, $x_{1}x_{4} =
x_{4}x_{2}$, $x_{2}x_{3} = x_{3}x_{1}$,
$x_{3}x_{2} = x_{4}x_{1}$. We know that $T$ is a
monoid of I-type (see \cite[Section 4]{jes}) and
thus, by Theorem~\ref{elem}, $T =
\{(a,\phi(a))\mid a \in \FaM_{n}\}$ for some map
$\phi : \FaM_{n}\rightarrow \Sym_{n}$. Put
$\FaM_{4} = \langle
u_{1},u_{2},u_{3},u_{4}\rangle$, $\phi(u_{i}) =
\sigma_{i}$ and $x_{i} = (u_{i},\sigma_{i})$. The
defining relations allow us to discover the
action on $\FaM_{4}$. For example, $x_{1}x_{2} =
(u_{1},\sigma_{1})(u_{2},\sigma_{2}) =
(u_{1}\sigma_{1}(u_{2}),\sigma_{1}\sigma_{2})$
and $x_{3}x_{3} =
(u_{3},\sigma_{3})(u_{3},\sigma_{3}) =
(u_{3}\sigma_{3}(u_{3}),\sigma_{3}\sigma_{3})$.
Since $x_{1}x_{2} = x_{3}x_{3}$ we get that
$\sigma_{1}(u_{2}) = u_{3}$ and
$\sigma_{3}(u_{3}) = u_{1}$. Going through all
the defining relations we obtain that
$$\sigma_{1} = (23),\; \sigma_{2} = (14),\;
\sigma_{3} = (1243),\;  \sigma_{4} = (1342).$$
Clearly, $G = \{\phi(a)\mid a \in \FaM_{4}\}
\cong D_{8}$, the dihedral group of order $8$.
Let $B = \langle
u_{i}^{3},u_{i}^{2}u_{j},u_{i}u_{j}u_{k}\mid 1
\leq i \neq j \neq k \leq 4\rangle$. Then $B$ has
a group of quotients $BB^{-1} = \gr (u_{1}^{3},
u_{1}u_{4}^{-1}, u_{2}u_{4}^{-1},
u_{3}u_{4}^{-1}) = \{ a\in \Fa_{4}\mid
\left|a\right|\in 3\Z \}$, where $ \left| a
\right|$ denotes the natural (total) degree of
$a$. Clearly, $B$ is $G$-invariant. Hence, $S =
\{(b,\phi(b))\mid b\in B\}\subseteq  B \rtimes G$
is a monoid of IG-type. Note that $G$ is now
considered as a subgroup of $Aut(B)$. Because
$BB^{-1}$ has torsion-free rank $4$ and since $B$
(and thus $S$) has $20$ indecomposable elements,
it follows that $S$ is not of I-type.

We now give an example of a monoid $S$ of IG-type
so that $SS^{-1}$ has non-trivial periodic
elements. On the other hand, $SS^{-1}$ does not
contain non-trivial finite normal subgroups and
thus $K[SS^{-1}]$ (and $K[S]$) are prime algebras
(see for example \cite{pas,boekok}).

\begin{example} \label{torsionex}
Let $A=\langle u_{1},u_{2},u_{3},u_{4}\mid
u_{1}u_{2}=u_{3}u_{4}\rangle$ be the maximal
order as in Example~\ref{and}. Let  $D_{8} =
\langle  a,b\mid a^{4} = 1, b^{2} = 1, a^{3}b =
ba \rangle$, with $a = (1324)$ and $b= (12)$, the
Dihedral group of order $8$. So $D_{8}$ acts
naturally on $A$. In the semidirect product
$A\rtimes D_{8}$ consider the elements $x_{i} =
(u_{i},\sigma_{i})$, where $\sigma_{1}= (1324)$,
$\sigma_{2} = (12)$, $\sigma_{3} =(1423)$ and
$\sigma_{4} = (34)$ and let $S=\langle x_{1},
x_{2}, x_{3}, x_{4}\rangle$. Then $S$ is a monoid
of IG-type. Furthermore,  $SS^{-1}$ has
non-trivial periodic elements but $SS^{-1}$ does
not have non-trivial finite normal subgroups. So
$K[S]$ is a prime ring.
\end{example}

\begin{proof}
It is easily verified that $S = \{(a,\phi(a))\mid
a\in A\}$ and thus $S$ is a monoid of IG-type,
with $G=\{ \phi (a) \mid a\in A\} =D_{8}$. So,
$SS^{-1}\subseteq AA^{-1}\rtimes D_{8}$. Clearly,
$$(u_{3},\sigma_{3})(\sigma_{1}^{-1}(u_{1}^{-1}),\sigma_{1}^{-1})
= (u_{3}u_{2}^{-1},(12)(34)),$$ and, as $u_{4} =
u_{1}u_{2}u_{3}^{-1}$ in $AA^{-1}$, we have that
$(u_{3}u_{2}^{-1},(12)(34))^{2} = 1$. So
$SS^{-1}$ has non-trivial periodic elements.

We claim now that $SS^{-1}$ does not contain a
finite normal subgroup, or equivalently,
$K[SS^{-1}]$ is prime. Indeed, since $AA^{-1}$ is
torsion-free, it is readily verified that  finite
normal subgroups $N$ of $SS^{-1}$ must be such
that their natural projection onto $G=D_{8}$ are
contained in $\gr (a)$. Furthermore, it then
follows that $N$ contains a finite normal
subgroup of $G$ that is of order $2$. So $N$
contains a central element of order $2$. But
central elements in $SS^{-1}$ are of the form
$(u_{1}^{i}u_{2}^{i}, 1)$, so they are not
periodic. This proves the claim.
\end{proof}

We finish this section by showing that the
infinite dihedral group $D_{\infty}$ is a group
of IG-type. It also gives an example with
non-trivial torsion.

\begin{example}
Let $\sigma$ be the non-trivial isomorphism of
the infinite cyclic group $\Z$. So $\sigma
(1)=-1$. Then $H = \{(a,\phi(a))\mid a\in \Z,\;
\phi(a) = \sigma \text{ if } a\in 2\Z + 1,\;
\phi(a) = 1 \text{ if } a\in 2\Z \}$ is a group
of IG-type and $H \cong D_{\infty}$.
\end{example}
\begin{proof}
Put $a = (2,1)$ and $b = (1,\sigma)$. Then, $H =
\gr (a,b)$ and as $bab^{-1} = a^{-1}, b^{2} = 1$
we have that $H \cong D_{\infty}$.
\end{proof}

In \cite{brown} it is shown that the group
algebra $K[D_{\infty}]$ is not a maximal order
and that this algebra is the key  in
characterizing when a group algebra of a
polycyclic-by-finite group is a prime maximal
order.

\section{Prime ideals and Maximal Orders}

Throughout this section $S = \{(a,\phi(a))\mid a
\in A\}\subseteq A\rtimes G$ is a monoid of
IG-type, with $A$ a finitely generated abelian
cancellative monoid, $G=\{ \phi (a) \mid a\in
A\}$ a finite group and $SS^{-1}$ is
torsion-free. For an ideal $I$ of $A$, we put
$(I,\phi(I)) = \{(a,\phi(a))\mid a \in I\}$. Note
that this is a right ideal of $S$. By $Spec(S)$
we denote the set of all prime ideals of $S$.
Recall that the height of $Q\in Spec(S)$ is, by
definition, the largest non-negative integer $n$,
so that $S$ has a chain of primes $Q_{0}\subset
Q_{1}\subset ... \subset Q_{n} = Q$. We denote
this height by $\hth (Q)$.

We first describe the prime ideals of $S$. For
this we will make use of the next theorem
(Theorem~1.4 in \cite{okni}). It is worth
mentioning (as is already done in \cite{okni})
that, since $SS^{-1}$ is a localization of $S$
with respect to an Ore set of regular elements of
Noetherian ring $K[S]$,  the prime ideals of the
group algebra $K[SS^{-1}]$ are in a one-to-one
correspondence with the prime ideals $P$ of
$K[S]$ that do not intersect $S$ (see for example
\cite[Theorem 9.22, Theorem 9.20 and Lemma
9.21]{warfield}). Since prime ideals of group
algebras of polycyclic-by-finite groups have been
well studied through the work of Rosablade (see
\cite{pas,pass}) we thus get a lot of information
on all prime ideals of $K[S]$.

\begin{proposition}                \label{min-primes}
Let $S$ be a submonoid of a  torsion-free
abelian-by-finite  group $G$ and let $K$ be a
field.
\begin{enumerate}
\item If $P$ is a prime ideal in $S$, then $K[P]$ is a prime ideal
      in $K[S]$.
\item If $Q$ is a prime ideal in $K[S]$ with $Q\cap S \neq \emptyset$,
      then $K[Q\cap S]$ is a prime ideal in $K[S]$.
\item The height one
      prime ideals of $K[S]$ intersecting $S$ are of the form
      $K[P]$, where $P$ is a minimal prime ideal of $S$.
\end{enumerate}
\end{proposition}

\begin{theorem} \label{primes-descrip}
Let $S = \{(a,\phi(a))\mid a \in A\}\subseteq
A\rtimes G$ be a monoid of IG-type and suppose
$SS^{-1}$ is torsion-free. The prime ideals $P$
of $S$ of height $m$ are the sets $(Q_{1}\cap ...
\cap Q_{n},\phi(Q_{1}\cap ... \cap Q_{n}))$ so
that
\begin{enumerate}
\item  each $Q_{i}$ is a prime ideal of $A$
of height $m$,
\item  $a\phi(a)(Q_{1}\cap ... \cap
Q_{n})\subseteq Q_{1}\cap ... \cap Q_{n},$ for
every $a\in A$(that is, $(Q_{1}\cap ... \cap
Q_{n},\phi(Q_{1}\cap ... \cap Q_{n}))$ is an
ideal of $S$),
\item condition (2) is not satisfied for an intersection
over a proper subset of $\{Q_{1},...,Q_{n}\}$
(that is, $(Q_{1}\cap ... \cap
Q_{n},\phi(Q_{1}\cap ... \cap Q_{n}))$ is a
maximal set satisfying conditions (1) and (2)).
\end{enumerate}
\end{theorem}
\begin{proof}
Let $P$ be a prime ideal of $S$ and let $K$ be a
field. Because of Proposition~\ref{min-primes},
$K[P]$ is a prime ideal of $K[S]$. Let $A^{k} =
\{a^{k}\mid a\in A\}$, where $k$ is a divisor of
the order of the group $G$ such that $\phi(a^{k})
= 1$, for every $a\in A$ (see
Lemma~\ref{laatste}). We identify the group
$A^{k}A^{-k}$ with its natural image in
$SS^{-1}$. The algebra $K[S]$ has a natural
gradation by the finite group
$SS^{-1}/A^{k}A^{-k}$. The homogeneous component
of degree $e$ (the identity of $A^{k}A^{-k}$) is
the semigroup algebra $K[A^{k}]$. So, by
Theorem~17.9 in \cite{pass}, $$K[P]\cap K[A^{k}]
= P_{1}\cap...\cap P_{n},$$ an intersection of
primes $P_{i}$ of $K[A^{k}]$, each of the same
height as $K[P]$ (these are all the primes of
$K[A^{k}]$ minimal over $K[P\cap A^{k}]$).
Clearly, $$P\cap A^{k} =
\bigcap_{i=1}^{n}(P_{i}\cap A^{k}),$$ each
$P_{i}\cap A^{k}$ is a prime ideal of $A^{k}$ and
thus $K[P_{i}\cap A^{k}]$ is a prime ideal of
$K[A^{k}]$. From Proposition~\ref{min-primes} and
Theorem~17.9 in \cite{pass} we verified that
$\hth (P) = \hth (P_{i}\cap A^{k})$. As every
$P_{i}$ is minimal over $K[P\cap A^{k}]$ and
because $K[P_{i}\cap A^{k}]$ also is a prime over
$K[P\cap A^{k}]$ and it is contained in $P_{i}$
it follows that $K[P_{i}\cap A^{k}] = P_{i}$. So
$$P\cap A^{k} =Q_{1}^{(k)}\cap \cdots \cap
Q_{n}^{(k)},$$ with $Q_{i}^{(k)} = P_{i}\cap
A^{k}$.

We also make another remark. Let $Q$ be a prime
ideal of $A$, then $Q^{(k)} = \{q^{k}\mid q\in
Q\}\subseteq Q\cap A^{k}$ and $Q^{(k)}$ is a
prime ideal of $A^{k}$. Furthermore $Q\cap A^{k}$
is a nil ideal modulo $Q^{(k)}$. Since $A^{k}$ is
commutative it follows that $Q\cap A^{k}
\subseteq Q^{(k)}$. Hence $Q\cap A^{k} =
Q^{(k)}$. So we have a bijection between the
primes of $Q$ and $Q^{(k)}$ (and corresponding
primes have the same height).

If $(a,\phi(a))\in P$ then
$$(a,\phi(a))(\phi(a)^{-1}(a),\phi(\phi(a)^{-1}(a)))...
= (a^{k},1)\in P\cap A^{k}.$$ Hence $a^{k}\in
Q_{1}^{(k)}\cap \cdots \cap Q_{n}^{(k)}$, and
thus $a\in Q_{1}\cap \cdots \cap Q_{n}$.
Therefore $P\subseteq (Q_{1}\cap \cdots \cap
Q_{n},\phi(Q_{1}\cap \cdots \cap Q_{n}))$.
Conversely, if $(b,\phi(b))\in (Q_{1}\cap \cdots
\cap Q_{n},\phi(Q_{1}\cap \cdots \cap Q_{n}))$,
then $$((b,\phi(b))^{k})^{k}\in (Q_{1}\cap \cdots
\cap Q_{n},\phi(Q_{1}\cap \cdots \cap Q_{n}))\cap
A^{k}\subseteq (\cap_{i=1}^{n} (Q_{i}\cap
A^{k}),1)\subseteq P.$$ So $(Q_{1}\cap \cdots
\cap Q_{n},\phi(Q_{1}\cap \cdots \cap Q_{n}))$ is
a right ideal of $S$ that is nil modulo $P$.
Since $S/P$ satisfies the ascending chain
condition on one sided ideals, it follows that
$(Q_{1}\cap \cdots \cap Q_{n},\phi(Q_{1}\cap
\cdots \cap Q_{n}))\subseteq P$ (see for example
17.22 in \cite{faith}). Hence $P = (Q_{1}\cap
\cdots \cap Q_{n},\phi(Q_{1}\cap \cdots \cap
Q_{n}))$.

Since $P$ is a left ideal we also have that
$a\phi(a)(Q_{1}\cap ... \cap Q_{n})\subseteq
Q_{1}\cap ... \cap Q_{n}$.

Next we show that if $P_{1} = (Q_{1}\cap \cdots
\cap Q_{n},\phi(Q_{1}\cap \cdots  \cap Q_{n}))$
and $P_{2} = (Q' _{1}\cap \cdots  \cap Q'
_{m},\phi(Q' _{1}\cap \cdots  \cap Q' _{m}))$ are
different prime ideals (of the same height) of
$S$ then $\{Q_{1},\cdots ,Q_{n}\}\cap \{Q' _{1},
\cdots  , Q' _{m}\} = \emptyset$. Indeed, suppose
the contrary, then, without loss of generality,
we may assume that $Q_{1} = Q' _{1}$. As
$P_{1}\neq P_{2}$, and because they are of the
same height, we thus get that say $n\rangle1$ and
$m\rangle1$.

Clearly  $(Q_{2}\cap \cdots  \cap
Q_{n},\phi(Q_{2}\cap \cdots \cap Q_{n}))$ is a
right ideal of $S$ and
 \begin{eqnarray*}
\lefteqn{ (Q_{2}\cap \cdots  \cap
Q_{n},\phi(Q_{2}\cap \cdots \cap Q_{n}))P_{2}}\\
&\subseteq& \{(a\phi(a)(Q' _{1}\cap \cdots \cap
Q'_{m}),\phi(a\phi(a)(Q'_{1}\cap \cdots  \cap Q'
_{m})))\mid a\in Q_{2}\cap \cdots \cap Q_{n}\}\\
&\subseteq&  ((Q' _{1}\cap \cdots  \cap Q'_{m})
\cap (Q_{2}\cap \cdots \cap Q_{n}),\phi((Q'
_{1}\cap \cdots  \cap Q' _{m})\cap (Q_{2}\cap
\cdots  \cap Q_{n})))\\ &\subseteq& P_{1}
 \end{eqnarray*}
But as $\hth (Q_{1}) = \hth (Q_{2})= ... = \hth
(Q_{n})$ and the primes $Q_{1},...,Q_{n}$ are
distinct it follows that $$Q_{2}\cap ... \cap
Q_{n} \nsubseteq Q_{1}\cap ... \cap Q_{n}.$$ As
$P_{1}$ is prime we thus get that $P_{2}\subseteq
P_{1}$. But since they are of the same height, it
follows that $P_{1} = P_{2}$, a contradiction.
The above claim of course implies the minimality
as stated in the Theorem.

To end the proof, we need to show that ideals
$(Q_{1}\cap \cdots  \cap Q_{n},\phi(Q_{1}\cap
\cdots \cap Q_{n}))$ with the listed properties
are prime ideals of $S$. We know that
$Q_{1}^{(k)} = Q_{1}\cap A^{k}$ is a prime ideal
of $A^{k}$ of the same height as $Q_{1}$. Also
$K[Q_{1}^{(k)}]$ is a prime ideal of $K[A^{k}]$.
Again using graded techniques and
\cite[17.9]{pass} we know that there exists a
prime ideal $P$ of $K[S]$ that lies over
$K[Q_{1}^{k}]$. So $K[Q_{1}^{k}]$ is a minimal
prime over $P\cap K[A^{k}]$ and $P\cap K[A^{k}]=
K[Q_{1}^{k}]\cap X_{2}\cap \cdots \cap X_{m}$,
where $X_{2},\cdots ,X_{m}$ are minimal primes
over $P\cap K[A^{k}]$ and they are of the same
height as $K[Q_{1}^{k}]$. Clearly $P\cap S$ is a
prime ideal of $S$ and $(P\cap S)\cap A^{k} =
Q_{1}^{(k)}\cap (X_{2}\cap A^{k})\cap \cdots \cap
(X_{m}\cap A^{k})$. Hence, by the first part of
the proof, $P_{1} = P\cap S = (Q_{1}\cap Q'
_{2}\cap \cdots \cap Q' _{m},\phi(Q_{1}\cap Q'
_{2}\cap ...\cap Q' _{m}))$ with $(Q_{i}^{(k)})'
= (X_{2} \cap A^{k})$.

We can now do the same for $Q_{2},\ldots ,Q_{n}$.
Hence, we get primes $P_{2},\ldots ,P_{n}$ of $S$
so that $P_{2} = (Q_{2}\cap J_{2},\phi(Q_{2}\cap
J_{2}))$,\ldots , $P_{n} = (Q_{n}\cap
J_{n},\phi(Q_{n}\cap J_{n}))$, with each $J_{i}$
an intersection of primes of $A$ that are of the
same height as $Q_{i}$.  Furthermore, because of
the assumptions we get that
$$(J_{2},\phi(J_{2}))(Q_{1}\cap...\cap
Q_{n},\phi(Q_{1}\cap...\cap Q_{n}))\subseteq
(J_{2}\cap Q_{1}\cap ...\cap Q_{n})\subseteq
P_{2}.$$ Since $(J_{2},\phi(J_{2}))\nsubseteq
P_{2}$ this yields that $(Q_{1}\cap...\cap
Q_{n},\phi(Q_{1}\cap...\cap Q_{n}))\subseteq
P_{2} = (Q_{2}\cap J_{2},\phi(Q_{2}\cap J_{2}))$.
Since $\{ Q_{1}, \ldots , Q_{n}\}$ satisfies the
minimality condition as stated in the Theorem, we
obtain that $Q_{2}\cap  J_{2} = Q_{1}\cap...\cap
Q_{n}$. Thus $P_{1} = P_{2} = .... = P_{n} =
(Q_{1}\cap ... \cap Q_{n},\phi(Q_{1}\cap ... \cap
Q_{n}))$ and thus this is a prime ideal of $S$.
\end{proof}

The following is an immediate consequence from
the previous result (and its proof).

\begin{corollary}
If $L = \{Q_{1},...,Q_{n}\}$ is a full $G$-orbit
of primes of the same height in $A$, then there
exists a partition $\{X_{1},...,X_{n}\}$ of $L$,
so that $$(\cap_{Q\in X_{i}} Q,\phi(\cap_{Q\in
X_{i}} Q)) = P_{i}$$ are prime ideals of $S$.
\end{corollary}

We now can prove the main result. It provides a
characterization of  semigroup algebras $K[S]$ of
monoids of IG-type that are a maximal order.

\begin{theorem}\label{hoofd}
Let $S = \{(a,\phi(a))\mid a \in A\}\subseteq
A\rtimes G$ be a monoid of IG-type. Suppose that
$SS^{-1}$ is torsion-free and suppose that the
abelian monoid  $A$ is finitely generated  and a
maximal order. Then, the Noetherian PI-domain
$K[S]$ is a maximal order if and only if the
minimal primes of $S$ are of the form $$P =
(Q_{1}\cap \cdots \cap Q_{n},\phi(Q_{1}\cap
\cdots \cap Q_{n})),$$ where $\{Q_{1}.... Q_{n}\}
= \{\phi(a)(Q_{1})\mid a\in A\}\subseteq
\Spec^{0}(A)$.
\end{theorem}

\begin{proof} Because of the assumption,
the results in the first section show that $K[S]$
is a Noetherian domain that satisfies a
polynomial identity. In particular, $S$ satisfies
the ascending chain condition on one sided
ideals.

We first prove the sufficiency of the mentioned
condition. So, suppose that the minimal primes of
$S$ are of the from $(Q_{1}\cap .... \cap
Q_{n},\phi(Q_{1}\cap .... \cap Q_{n})),$ where
$\{Q_{1}.... Q_{n}\} = \{\phi(a)(Q_{1})\mid a\in
A\}\subseteq \Spec^{0}(A)$. To prove that $K[S]$
is a maximal order, it is sufficient to verify
conditions (2) and (3) of Theorem~\ref{tool}. The
former says that $S$ is a maximal order. Because
of Lemma~4.4 in \cite{okni}, in order to prove
this property,  it is sufficient to show that
$(P:_{l}P)=(P:_{r}P)=S$ for every prime ideal $P$
of $S$. From Theorem~\ref{primes-descrip} we know
that $P=(Q,\phi (Q))$ with $Q$ an intersection of
prime ideals in $A$ of the same height, say $n$.
Assume $(x,\phi (x)) \in (P:_{l}P)$. Then  $x\phi
(x)(Q) \subseteq Q$. If $n\neq 0$ (so $Q$ is an
intersection of primes that are not minimal) then
the divisorial closure of both $Q$ and  $\phi
(x)(Q)$ equals $A$. As $x(\phi (x)(Q))^{*}
\subseteq Q^{*}$ we thus get that $x\in A$ and
thus $(x,\phi (x))\in S$. If $n=0$, then, by
assumption, $Q$ is $G$-invariant and thus we get
that $xQ\subseteq Q$. Since $A$ is a maximal
order, this yields that $x\in A$ and again
$(x,\phi (x))\in S$. So $(P:_{l}P)=S$. On the
other hand, suppose $(Q,\phi (Q)) (x,\phi
(x))\subseteq (Q,\phi (Q))$. Then $A(Q\cap
A^{k})x\subseteq Q$. If $n\neq 0$ then $Q\cap
A^{k}$ is not contained in a minimal prime ideal
of $A$ and thus the divisorial closure of both
$A(Q\cap A^{k})$ and $Q$ is $A$. Since $(A(Q\cap
A^{k}))^{*} x\subseteq Q^{*}$ we get that
$(x,\phi (x))\in S$. So it remains to show that
$(P:_{r}P)=S$ for a minimal prime ideal $P$.
Hence, by assumption $P =(Q,\phi (Q))$ with $Q$ a
$G$-invariant ideal. More generally, we prove
that $(I:_{r}I)=S$ for any ideal $I=(M,\phi (M))$
of $S$ with $M$ a $G$-invariant ideal of $A$.

We now prove this by contradiction. So suppose
that $I$ is such an ideal of $S$ with
$Ig\subseteq I$ for some $g\in SS^{-1}\setminus
S$.

Now, as $g\in SS^{-1}$, we know that $g =
(a,\phi(a)) (z,1)^{-1}$ with $z$ an invariant
element of $A$, and thus $(z,1)$ central in $S$.
As $A$ is a maximal order, we have that the
minimal primes of $A$ freely generate the abelian
group $D(A)$. So, in the divisor group $D(A)$, we
can write $Az$ as a product of minimal primes.
Because $Az$ is invariant, the minimal primes in
a $G$-orbit have the same exponent. Hence,
 $$Az = \left( J_{1}^{n_{1}} \right)^{*} * \cdots
  *  \left (J_{l}^{n_{l}}\right)^{*},$$
where each $J_{i}$ is an intersection of all
minimal primes of $A$ in a $G$-orbit. So,
 because of the
assumption and Theorem~\ref{primes-descrip}, each
$(J_{i},\phi(J_{i}))$ is a minimal prime of $S$.
Of course also $Aa$ is a divisorial product of
minimal primes of $A$. If necessary, cancelling
some common factors of $Aa$ and $Az$, we may
assume  that $Aaz^{-1}= K*L^{-1}\not\subseteq A$,
and thus $KL^{-1} \not\subseteq A$ with $L=\left(
J_{1}^{n_{1}}\right)^{*} * \cdots * \left(
J_{l}^{n_{l}}\right)^{*}$, $L^{-1} = (A:L)$ and
$K$ is not contained in $J_{i}$, for every $i$
with $1\leq i \leq l$. Note that, also, $L^{-1}$
is $G$-invariant and thus $(L^{-1},\phi(L^{-1}))$
is a fractional ideal of $S$. Of course,
$I(K,\phi (K))(L^{-1},\phi (L^{-1})\subseteq I $.
Because $S$ satisfies the ascending chain
condition on ideals, we can choose $I$ maximal
with respect to the property that such $K$ and
$L$ exist with $KL^{-1} \not\subseteq A$.

Clearly we obtain that
\begin{eqnarray*}
 I (K,\phi (K)) (L^{-1},\phi (L^{-1})) (L,\phi (L))
  &\subseteq&
 I (L,\phi (L)) \subseteq S (L,\phi (L))\\
 & \subseteq &
 (J_{i}^{n_{i}},\phi(J_{i}^{n_{i}}))\\
 & \subseteq&
 (J_{i},\phi(J_{i})).
\end{eqnarray*}
Since, $(J_{i},\phi(J_{i}))$ is a prime ideal of
$S$, we get that either $(KL^{-1}L,\phi
(KL^{-1}L)) \subseteq  (J_{i},\phi(J_{i}))$ or $I
\subseteq (J_{i},\phi(J_{i}))$. Because of the
above, the former is excluded. Hence $I \subseteq
(J_{i},\phi(J_{i}))$. As $J_{i}$ is
$G$-invariant, we get again that
$(J_{i}^{-1},\phi(J_{i}^{-1}))$ is a fractional
ideal of $S$ that contains $S$. Therefore, we get
that $I \subseteq (J_{i}^{-1},\phi(J_{i}^{-1}))
I$ is an ideal of $S$. Since
$$(J_{i}^{-1},\phi(J_{i}^{-1})) I  (K,\phi (K))
(L^{-1},\phi (L^{-1})) \subseteq
(J_{i}^{-1},\phi(J_{i}^{-1}))I,$$ the maximality
condition on $I$ thus implies that
$$(J_{i}^{-1},\phi(J_{i}^{-1}))I = I.$$ Since $M$
is $G$-invariant this yields that
 $J_{i}^{-1}M =
M$ and thus $J_{i}^{-1}*M^{*} = M^{*}$. So
$J_{i}^{-1} =A$, a contradiction.

We now show that, if $P$ is a minimal prime ideal
of $S$ then the monoid $S_{P}$ has only one
minimal prime. As $A^{k}$ is $G$-invariant,
Lemma~2.3 in \cite{okni} gives that $S_{P} =
S_{(P)}$ with
 \begin{eqnarray*}
  S_{(P)} &=& \{g\in SS^{-1}\mid Cg \subseteq S, \;
 C\not\subseteq P, \; C\subseteq A^{k},\\
 &&
\text{ for some conjugacy class } C \text{ of }
SS^{-1}\}.
 \end{eqnarray*}
Again, by assumption, $P = (Q_{1}\cap...\cap
Q_{n},\phi(Q_{1}\cap...\cap Q_{n}))$, where
$\{Q_{1},...,Q_{n}\}$ is a full $G$-orbit of
minimal primes of $A$. First we prove that
 \begin{eqnarray}  \label{form-loc}
  S_{(P)} &=&
(A_{Q_{1}}\cap...\cap
A_{Q_{n}},\phi(A_{Q_{1}}\cap...\cap A_{Q_{n}})) .
 \end{eqnarray}
So suppose that $(b,\phi(b))\in S_{(P)}$. Then
there exists a conjugacy class $C$ of $SS^{-1}$
in $A^{k}$ with $$C(b,\phi(b)) \subseteq S$$ and
$C$ not contained in $P$. But as $C$ is contained
in $A^{k}$ it is easily verified that $C =
\{(\phi(a)(c^{k}),1)\mid \phi(a)\in G\},$ for
some $c\in A$. Since $C \nsubseteq P$ and because
$\{Q_{1},...,Q_{n}\}$ is a $G$-orbit, this yields
that $C \nsubseteq Q_{i}$, for every $i \in
\{1,...,n\}$. Hence, it follows that $b \in
A_{Q_{1}}\cap...\cap A_{Q_{n}}$. Conversely,
suppose that $b\in A_{Q_{1}}\cap ... \cap
A_{Q_{n}}$. Then there exist $c_{i}\in A\setminus
Q_{i}$ with $c_{i}^{k}b\in A$ for every $i \in
\{1,\ldots ,n\}$. As $M = A^{k} c_{1}^{k}\cup
\cdots  \cup A^{k} c_{n}^{k} \nsubseteq
Q_{1}^{(k)}\cap \cdots \cap Q_{n}^{(k)}$, it
follows that, in $D(A^{k})$, $M^{*}$ is a product
of minimal primes that do not belong to
$\{Q_{1}^{(k)},\ldots,Q_{n}^{(k)}\}$. So, $N =
\prod_{g\in G} g(M^{*})$ is an invariant ideal of
$A^{k}$ and $N\not\subseteq Q_{1}\cap \cdots \cap
Q_{n}$. Clearly, $NbA \subseteq A$. Choose $d\in
N\setminus P$. Then $C' =\{\phi(a)(d)\mid a \in
A\}$ is a $SS^{-1}$- conjugacy class contained in
$A^{k}$, but not in $P$. Since, $C'
(b,\phi(b))\subseteq S$, we get that
$(b,\phi(b))\in S_{(P)}$, as desired. This
finishes the proof of (\ref{form-loc}).

The monoid $B=A_{Q_{1}}\cap \cdots \cap
A_{Q_{n}}$ is a maximal order with minimal prime
ideals $P_{i} = A_{Q_{1}}\cap ...\cap Q_{i}A_{Q_{i}}\cap... \cap A_{Q_{n}}$,
$1\leq i \leq n$. From Lemma~2.2
in \cite{okni} we know that $I(P)=\{ (x,\phi (x))
\in S_{P} \mid (x,\phi (x)) C \subseteq P, \text{
for some } G-\text{conjugacy class } C\subseteq S
\text{ with } C\not\subseteq P\}$ is a prime
ideal of $S_{P}$. It is easily seen that
 $I(P)=\{ (x,\phi (x)) \in
S_{P} \mid (x,\phi (x)) C \subseteq P, \text{ for
some } G-\text{conjugacy class } C\subseteq A^{k}
\text{ with } C\not\subseteq P\}$. From
(\ref{form-loc}) it then follows that
 $I(P)=(BQ_{1}\cap \cdots \cap BQ_{n}, \phi
 (BQ_{1}\cap \cdots \cap BQ_{n}))$.
Therefore, Theorem~\ref{primes-descrip} implies
that this is the only  minimal prime ideal of
$S_{P}$. This finishes the proof of the
sufficiency of the conditions.

To prove the necessity, assume $K[S]$ is a
maximal order.  Let $P = (M,\phi (M))$ be a
minimal prime ideal of $S$. Theorem~\ref{tool}
yields that $S_{P}$ has a unique minimal prime.
Furthermore, since $A^{k}$ is $G$-invariant,
Lemma~2.5 in \cite{okni} yields that $(M,\phi
(M))\cap A^{k}=M\cap A^{k}$ is $G$-invariant.
Consequently, $M$ is $G$-invariant and
Theorem~\ref{primes-descrip} yields that $M$ is
the intersection of a full $G$-orbit of minimal
primes. This finishes the proof.
\end{proof}

\begin{remark} \label{remark-torsion}
Let $S = \{(a,\phi(a))\mid a \in A\}\subseteq
A\rtimes G$ be a monoid of IG-type. Suppose that
the abelian monoid  $A$ is finitely generated and
a maximal order. If the minimal primes of the
monoid $S$ are  as stated in the sufficient
condition of Theorem~\ref{hoofd} then  $S$ is a
maximal order. (The assumption $SS^{-1}$ is
torsion-free, is not needed in the proof of this
part of the result.)
\end{remark}

As an application of Theorem~\ref{hoofd} we give
two examples.

\begin{example}
Let $S$ be the  monoid of $IG$-type considered in
Example~\ref{and}. The semigroup $S$ is a maximal
order and  $P_{1} =(Q_{1}\cap
Q_{4},\phi(Q_{1}\cap Q_{4}))$ and $P_{2}
=(Q_{2}\cap Q_{3},\phi(Q_{2}\cap Q_{3}))$ are its
minimal prime ideals. Furthermore,  $K[S]$ is a
maximal order for any field $K$.
\end{example}
\begin{proof} From
Example~\ref{and} (and its proof) we know that
$A$ is a finitely generated maximal order with
four minimal primes: $Q_{1} = (u_{1}, u_{3})$,
$Q_{2} = (u_{1}, u_{4})$, $Q_{3} = (u_{2},
u_{3})$ and $Q_{4} = (u_{2}, u_{4})$. Because
$SS^{-1}$ is torsion-free,
Proposition~\ref{primes-descrip} yields a
description of the prime ideals of $S$. Clearly,
$a\phi(a)(P_{i})\subseteq P_{i}$, for $i\in
\{1,2\}$. Hence because of Theorem~\ref{hoofd},
to  prove that $P_{1}$ and $P_{2}$ are the only
minimal primes of $S$, it is now sufficient to
note that for every $Q_{i}$, there exists an
$a\in A$ such that $a\phi(a)(Q_{i})\nsubseteq
Q_{i}$. Indeed, $u_{2}(12)(34)Q_{1}\nsubseteq
Q_{1},$ $u_{2}(12)(34)Q_{2}\nsubseteq Q_{2},$,
$u_{4}(12)(34)Q_{3}\nsubseteq Q_{3}$, and
$u_{3}(12)(34)Q_{4}\nsubseteq Q_{4}$.
\end{proof}

\begin{example}
Let $S$ be the monoid of IG-type defined in
Example~\ref{belvb}. Then $K[S]$ is a maximal
order for any field $K$.
\end{example}
\begin{proof}
It is readily verified that the minimal primes of
$S$ are of the form as required in
Theorem~\ref{hoofd}.
\end{proof}

We finish  this paper with an example of a monoid
of IG-type that is not a maximal order.

\begin{example}
Let $A=\langle u_1, u_2, u_3 , u_4 \mid u_1 u_2
u_3  = u_{4}^{2} \rangle$. Then $A$ is a maximal
order (in its torsion-free group of quotients)
with minimal prime ideals $Q_{1}=(u_{1},u_{4})$,
$Q_{2}=(u_2, u_4)$ and  $Q_{3} =(u_{3},u_{4})$.
Let $S=\{ (a,\phi (a)) \mid a\in A\}$, with $\phi
(a) =1$ if $a\in A$ has even degree in $u_4$,
otherwise, $\phi (a) = (1 2)$ (the transposition
interchanging $u_1$ with $u_2$). Then,
$S\subseteq A\rtimes \Z_{2}$ is a monoid of
IG-type which is not a maximal order and the
group of quotients $SS^{-1}$ is torsion-free.
Thus,  $K[S]$ is not a maximal order for any
field $K$. The minimal prime ideals of $S$ are
$P_{i}=(Q_{i},\phi (Q_{i}))$ with $1\leq i \leq
3$.
\end{example}

\begin{proof}
Clearly, $AA^{-1}=\gr (u_{1},u_{2},u_{4})$ is a
free abelian group of rank $4$. Furthermore,
every element of $AA^{-1}$ has a unique
presentation of the form
$u_{1}^{i}u_{2}^{j}u_{3}^{k}u_{4}^{m}$, with
$m\in \{ 0,1\}$ and $i,j,k \in \Z$; elements of
$A$ are those with $i,j,k$ non-negative. It is
easily seen that $Q_{1},Q_{2}$ and $Q_{3}$ are
the minimal primes of $A$ and the localizations
of $A$ with respect to these prime ideals are
$A_{Q_{1}}= A\langle
u_{2}^{-1},u_{3}^{-1}\rangle$, $A_{Q_{2}}=
A\langle u_{1}^{-1},u_{3}^{-1}\rangle$ and
$A_{Q_{3}}= A\langle
u_{1}^{-1},u_{2}^{-1}\rangle$. Furthermore,
$A_{Q_{1}}\cap A_{Q_{2}} \cap A_{Q_{3}} =A$ and
each $A_{Q_{i}}$ is a maximal order with unique
minimal prime ideal $u_{4}A_{Q_{i}}$. It follows
that $A$ is a maximal order. As $\Z_{2}=\gr
((12))$ induces a faithful action on the finitely
generated monoid $A=\langle u_1, u_2, u_3,
u_4\mid u_1 u_2 u_3 = u_{4}^{2} \rangle$, we thus
get that $S\subseteq A\rtimes \Z_{2}$ is a monoid
of IG-type.

Suppose that there exists a non-trivial  periodic
element in the group of quotients $SS^{-1}$. Such
an element must be of the form
$(u_{1}^{\alpha_{1}}u_{2}^{\alpha_{2}}u_{4}^{\alpha_{4}},
(12))$, with $\alpha_{i}\in \Z$ and $\alpha_{4}$
odd. Then, by Theorem~\ref{karel} (see also the
remarks stated after its proof),
$(\alpha_{1},\alpha_{2},\alpha_{4}) + (12)
(\alpha_{1}', \alpha_{2}', \alpha_{4}') =
(\alpha_{1}', \alpha_{2}', \alpha_{4}')$, for
some $\alpha_{1}', \alpha_{2}', \alpha_{4}'\in
\R$. Hence $\alpha_{4} + \alpha_{4}' =
\alpha_{4}'$ and thus $\alpha_{4}=0$, a
contradiction. So, $SS^{-1}$ indeed is
torsion-free.

Let $I$ be the ideal generated by $((u_{1},1),
(u_{4},(12)))$. Then $(u_{1}u_{3}u_{4}^{-1},(12))
I \subseteq I $ and thus $S$ is not a maximal
order.

As an immediate consequence of
Theorem~\ref{primes-descrip} one sees that $P_{1}
= (Q_{1},\phi (Q_{1}))$, $P_{2} = (Q_{2},\phi
(Q_{2}))$  and $P_{3}= (Q_{3},\phi (Q_{3}))$ are
the minimal prime ideals of $S$. Clearly $Q_{1}$
and $Q_{2}$ are not $\Z_{2}$-invariant, as
corresponds with Remark~\ref{remark-torsion}.
\end{proof}

 \noindent \\
 Department of mathematics\\
 Vrije Universiteit Brussel\\
 Pleinlaan 2\\
 1050 Brussel, Belgium\\
 email: efjesper@vub.ac.be and igoffa@vub.ac.be


\begin{thebibliography}{99}
\itemsep=-2pt
\bibitem{and} D.F. Anderson,
The divisor class group of a semigroup ring,
Comm. Algebra 8(5) (1980), 467-476.
\bibitem{bab} A.Babakhanian,
Cohomological Methods in Group Theory, Marcel
Dekker, New York, 1972.
\bibitem{brown} K.A. Brown,
Height one primes of polycyclic group rings,
London, Math. Soc, 1985. (426-438).
\bibitem{3} F. Cedo, E.Jespers and  J.Okni\'{n}ski,
Semiprime quadratic algebras of Gelfand-Kirillov
dimension one. J. Algebra Appl. 3 (2004),
283-300.
\bibitem{clifpres} A.H.Clifford and G.B.Preston,
The Algebraic Theory of Semigroups, Vol.I,
Amer.Math.SOC., Providence, 1961.
\bibitem{choui} II. L.G. Chouinard,
Krull semigroups and divisor class groups, Canad.
J. Math. 23 (1981), 1459-1468.
\bibitem{karel} K. Dekimpe, Private communication.
\bibitem{eting-gur} Etingof P.,  Guralnick R. and
Soloviev A., Indecomposable set-theoretical
solutions to the quantum Yang-Baxter equation on
a set with a prime number of elements, J. Algebra
249 (2001), 709--719.
\bibitem{eting}  Etingof P., Schedler T. and
Soloviev A., Set-theoretical solutions of the
quantum Yang-Baxter equation, Duke Math. J. 100
(1999), 169--209.
\bibitem{faith} C. Faith, Algebra II, Ring Theory,
Springer-Verlag, New York, 1976.
\bibitem{3b} T. Gateva-Ivanova, E.Jespers and J.
Okni\'{n}ski, Quadratic algebras of skew type and
the underlying monoids, J. Algebra 270 (2003),
635-659.
\bibitem{gateva} T. Gateva-Ivanova and M. Van den
Bergh, Semigroups of I-type, J.Algebra 206
(1998), no. 1, 97-112.
\bibitem{gilmer} Gilmer R., Commutative semigroup
rings, Univ.Chicago Press, Chicago, 1984.
\bibitem{warfield} Goodearl K.R. and Warfield R.B.,
An Introduction to Noncommutative Noetherian
Rings, Cambridge Univ. Press, New York, 1989.
\bibitem{noet} E.Jespers and  J. Okni\'{n}ski,
Noetherian semigroup algebras of submonoids of
polycyclic-by-finite groups, Bull. London Math.
Soc., to appear.
\bibitem{jes} E.Jespers and J.Okni\'{n}ski,
Monoids and groups of I-type, Algebras Repres.
Theory 8 (2005), 709-729.
\bibitem{okni} E.Jespers, J.Okni\'{n}ski,
Semigroup algebras and Noetherian maximal orders,
J.Algebra 238 (2001), 590-662.
\bibitem{pbf} E.Jespers and  J.Okni\'{n}ski,
Submonoids of Polycyclic-by-finite Groups and
their Algebras, Algebras and Representation
Theory 4 (2001), 133-153.
\bibitem{jeswang} E.Jespers, Wang Q., Height one
prime ideals in semigroup algebras satisfying a
polynomial identity, J. Algebra 248
(2002),118-131.
\bibitem{boekok} J.Okni\'{n}ski, Semigroup Algebras,
Marcel Dekker, New York, 1991.
\bibitem{pas} Passman D.S., The Algebraic
Structure of Group rings, New York, 1977.
\bibitem{pass} Passman D.S., Infinite Crossed
Products, Academic Press, New York, 1989.
\end{thebibliography}
\end{document}